\documentclass[11pt,leqno]{article} 
\usepackage{graphics}

\newtheorem{lma}{Lemma}[section]

\newcommand{\beqa}{\begin{eqnarray}}
\newcommand{\eeqa}{\end{eqnarray}}

\newcommand{\pf}{\noindent {\bf Proof:} $\s$ }
\newcommand{\epf}{ \hfill$\diamondsuit$ \medskip}

\newcommand{\beq}{\begin{equation}}
\newcommand{\eeq}{\end{equation}}
\newcommand{\lbl}{\label}
\newcommand{\s}{\; \;}

\newcommand{\la}{\lambda}
\newcommand{\mb}{\mbox}
\newcommand{\ra}{\rightarrow}
\newcommand{\al}{\alpha}
\newcommand{\p}{\varphi}

\title{Continuation of global solution curves using global parameters }

\author{
Philip Korman   \\ 
Department of Mathematical Sciences \\ 
University of Cincinnati \\ 
Cincinnati Ohio 45221-0025 \\
\\
Dieter S. Schmidt \\
Department of Computer Science \\
University of Cincinnati \\
Cincinnati, Ohio 45221-0030}

\date{}

\begin{document}

\maketitle
\begin{abstract} 
This paper provides both the theoretical results and numerical calculations of global solution curves, by continuation in  global parameters. Each point on the solution curves is computed directly as the global parameter is varied, so that all of the turns that the solution curves make, as well as its different branches, appear automatically on the computer screen. For radial $p$-Laplace equations we present a simplified derivation of the regularizing transformation from P. Korman \cite{K2015}, and use this transformation for more accurate numerical computations. While for $p>2$ the solutions are not of class $C^2$, we show that they are of the form $w(r^{\frac{p}{2(p-1)}})$, where $w(z)$ is of class $C^2$. Bifurcation diagrams are also calculated for non-autonomous problems, and for the fourth order equations modeling elastic beams. We show that the first harmonic of the solution can also serve as a global parameter. 
 \end{abstract}

\begin{flushleft}
\end{flushleft}

\begin{flushleft}
\end{flushleft}

\section{Introduction}
\setcounter{equation}{0}
\setcounter{thm}{0}
\setcounter{lma}{0}

\medskip

We wish to compute the global bifurcation diagrams for the operator equation
\[
Lu=\la f(u) \,.
\]
Here $L$ is a linear operator, $\la $ is a parameter. A traditional approach involves continuation in $\la$, where one uses Newton's method, with the initial iterate being the solution at the preceding value of $\la$ (or its linear extrapolation to the present value of $\la$). These ``curve following" methods are very well developed, see e.g., the book of  E.L. Allgower  and K. Georg \cite{A}, but not easy to implement. The problem is that the solution set may consist of multiple pieces, and there may be multiple turns on each curve. If a turn occurs, your computer program must recognize that fact, jump on the other branch, and reverse the direction of continuation in $\la$. Here $\la $  is a {\em local parameter}. If, on the other hand, the problem possesses a {\em global parameter} uniquely identifying the solution pair
$(\la,u)$, then after computing sufficiently many of these pairs directly, one gets the picture of the entire solution set, including all of the turns on each piece.
\medskip

What quantity can serve as a global parameter? In this paper we present two possibilities: it is either the maximum value of the solution, or the first harmonic of the solution. Next we discuss the first of these possibilities for problems on a ball in $R^n$. 
\medskip

By  classical theorem of B. Gidas, W.-M. Ni and L. Nirenberg \cite{GNN}, any positive solution of the semilinear Dirichlet problem (here $u=u(x)$, $x \in R^n$, $\la$ is a positive parameter)
\[
\Delta u+ \la \, f(u) =0 \s \mb{for $|x|<1$}\,, \s u=0  \s \mb{when $|x|=1$}
\]
is necessarily radially symmetric, i.e., $u=u(r)$, with $r=|x|$, and so the problem 
turns into an ODE
\beq
\lbl{2.15.1a}
\s\s\s u''(r)+\frac{n-1}{r} u'(r) +\la f(u(r))=0 \s \mbox{for $0<r<1$}, \s u'(0)=u(1)=0 \,.	
\eeq
Moreover, this theorem asserts that the function $u(r)$ is strictly decreasing, and hence  $u(0)$ gives the maximum value of solution. A simple scaling argument shows that the value of $u(0)>0$ is a global parameter, uniquely identifying the solution pair
$(\la,u(r))$. We continue the solutions in $u(0)$ by developing the natural ``shoot-and-scale" method. This method is very easy to implement, even when the solution set of (\ref{2.15.1a}) consists of multiple curves, with multiple turns on each curve. This method  was used previously by R. Schaaf \cite{S10} for the one-dimensional case, $n=1$.  We also develop the ``shoot-and-scale" method in case of the Neumann boundary conditions. We include very short {\em Mathematica} programs, producing the global solution curves for both the Dirichlet and the Neumann boundary conditions. Similar computations were used in P.  Korman, Y. Li  and  D.S. Schmidt \cite{KLS} to produce  nodal ground state solutions of Schr\"odinger's equation.
\medskip

For non-autonomous problems
\[
u''(r)+\frac{n-1}{r} u'(r) +\la f(r,u(r))=0 \s\s \mbox{for $0<r<1$}, \s u'(0)=u(1)=0 	
\]
the ``shoot-and-scale" method does not apply (one can still shoot, but not scale). In particular, one cannot assert anymore that the value of $u(0)$ is a global parameter. It is plausible that in many cases the value of $u(0)$ is still a global parameter (we mention below such a case  when $n=1$). Under the assumption that $u(0)$ is a global parameter we develop a continuation method in $u(0)$, based on  Newton's method.
\medskip

We then develop  the shoot-and-scale method to compute the curves of positive solutions for the Dirichlet problem  in case of the $p$-Laplacian
\beq
\lbl{sept18}
\p (u')'+\frac{n-1}{r} \, \p (u') +\lambda f(u)=0, \s  u'(0)=u(1)=0 \,,
\eeq
where $\p (v)=v |v|^{p-2}$, $p>1$. As in case $p=2$, we show that the value of $u(0)$ is a global parameter, uniquely identifying the solution pair $(\lambda,u(r))$, so that the global  solution curves can be drawn in  the $(\lambda, u(0))$ plane. We perform numerical computations, and the corresponding {\em Mathematica} program is available from the authors. Writing the equation (\ref{sept18}) in the form 
\[
u''+\frac{n-1}{(p-1)r}u'+\lambda \frac{f(u)}{(p-1) |u'|^{p-2}}=0 \,,
\]
one sees that if  $p>2$, then $u''(0)$ does not exist, so that  $u(r) \not\in C^2[0,1)$. This singularity  also presents a difficulty for computations. Using a {\em regularizing transformation} 
\[
z=r^{\frac{p}{2(p-1)}} \,,
\]
it was shown in P. Korman \cite{K4} that $u=w \left(r^{\frac{p}{2(p-1)}} \right)$, where the function $w(z)$ is twice differentiable at $z=0$, $w(z) \in C^2[0,1)$. We present a simplified derivation of this result, and then use the regularizing transformation  for more efficient computation of the global solution curves.
\medskip

We  also develop  the shoot-and-scale method to compute the curves of positive solutions for a class of fourth order equations ($\la$ is a positive parameter)
\begin{eqnarray} 
\lbl{bm1a}
& u''''(x)=\la f(u(x)), \s  \mbox{for $x \in (-1,1)$} \\ \nonumber
&  u(-1)=u'(-1)=u(1)=u'(1)=0 \,, \nonumber
\end{eqnarray}
modeling the displacements of an elastic beam,  clamped at the end points. Under the assumption that $f(u)>0$ and $f'(u)>0$ for $u>0$, it was shown in P. Korman \cite{K} that  any positive solution is an even function, decreasing on $(0,1)$, so that $u(0)$ is the maximal value of the solution. Quite surprisingly, it was also shown in P. Korman \cite{K} that   the value of $u(0)$ is a global parameter, uniquely identifying the solution pair $(\la, u(x))$ of (\ref{bm1a}). This opens the way for continuation in $u(0)$, but the shoot-and-scale method requires a major adjustment. The reason is that the knowledge of $u(0)=\al $ and $u'(0)=0$ does not identify the solutions of the initial value problem corresponding to the equation in (\ref{bm1a}). We develop a method to compute the bifurcation diagrams for the problem (\ref{bm1a}).  Such curves have never  been computed previously, to the best of our knowledge.
\medskip

What other quantities can be used as a global parameter? Consider the Dirichlet problem
\beq
\lbl{sept18a}
\Delta u +f(u)=g(x) \,, \s x \in D \,, \s\s u(x)=0 \s \mbox{on $\partial D$} \,,
\eeq
on a smooth domain $D \subset R^n$. Decompose the given function  $g(x)=\mu \p _1(x)+e(x)$, where $\p _1(x)$ is the principal eigenfunction of the Laplacian with zero boundary conditions, $\mu \in R$, and $\int _D e(x) \p _1(x) \, dx=0$. Similarly, decompose the solution $u(x)=\xi \p _1(x)+U(x)$, with $\xi \in R$ and $\int _D U(x) \p _1(x) \, dx=0$. It was shown in P. Korman \cite{K2014} that $\xi$ is a global parameter, uniquely identifying the solution pair $(u(x), \mu)$ of (\ref{sept18a}), provided that 
\beq
\lbl{int2}
f'(u)<\la _2 \,, \s \mbox{ for all $u \in R$} \,,
\eeq
where $\la _2$ is the second eigenvalue of the Laplacian on $D$, with zero boundary conditions. We implement numerical computations of the solution curve $\left(u(x),\mu \right)(\xi)$ for the one-dimensional case, and apply them to illustrate some recent multiplicity results of A. Castro et al \cite{C}. Our computations also show that the condition (\ref{int2}) cannot be removed, in general.
\medskip

We see absolutely no need to ever use the traditional method of finite differences (or finite elements) for the problem (\ref{2.15.1a}). To explain why, it suffices to consider the one-dimensional version of the problem (\ref{2.15.1a}):
\begin{equation}
\label{2.15.3}
u''(x)+\la f(u(x))=0 \s\s \mb{for $-1<x<1$}, \s u(-1)=u(1)=0 \,.
\end{equation}
If we divide the interval $(-1,1)$ into $N$ pieces, with the step $h=2/N$, and the subdivision points $x_i=ih$, and denote by $u_i$ the numerical approximation of $u(x_i)$, then the finite difference approximation of (\ref{2.15.3}) is 
\begin{equation}
\label{2.15.4a}
\s\s \frac{u_{i+1}-2u_i+u_{i-1}}{h^2}+\la f(u_i)=0, \s 1 \leq i \leq N-1, \s u_0=u_N=0.
\end{equation}
This is a system of {\em nonlinear algebraic equations}, more complicated in every way than the original problem (\ref{2.15.3}). In particular, this system often has more solutions than the corresponding differential equation (\ref{2.15.3}). The existence of the extra solutions of difference equations (not corresponding to the solutions of the differential equation (\ref{2.15.3})) has been recognized for a while, and a term {\em spurious solutions} has been used. When studying the algebraic nonlinear system of equations  (\ref{2.15.4a}), one can no longer rely on the familiar tools from differential equations. Even for the simplest nonlinearity $f(u)=u^k$, the analysis of the finite difference problem (\ref{2.15.4a}) is very involved,  see E.L. Allgower \cite{Al0}. We avoid spurious solutions with the  methods used in this paper.

\section{Preliminary results}
\setcounter{equation}{0}
\setcounter{thm}{0}
\setcounter{lma}{0}

The study of radial solutions  begins with the classical theorem of B. Gidas, W.-M. Ni and L. Nirenberg \cite{GNN}, which states that in case of  a ball in $R^n$, say a unit ball,  any positive solution of the semilinear Dirichlet problem 
\begin{equation}
\label{2.1.0}
\Delta u+ \la  f(u) =0 \s \mb{for $|x|<1$}\,, \s u=0  \s \mb{when $|x|=1$}
\end{equation}
is necessarily radially symmetric, i.e., $u=u(r)$, with $r=|x|$, and so the problem 
turns into an ODE
\begin{equation}
\label{2.1.1}
\s\s\s u''(r)+\frac{n-1}{r} u'(r) +\la f(u(r))=0 \s \mbox{for $0<r<1$}, \s u'(0)=u(1)=0 \,.	
\end{equation}
Moreover, this theorem asserts that
\begin{equation}
\label{2.1.2}
u'(r)<0, \s\s \mbox{for all $0<r<1$} \,.
\end{equation}
It follows that $u(0)$ gives the maximum value of the solution $u(r)$.	
The only assumption of this remarkable theorem is a slight smoothness assumption on $f(u)$, which is considerably weaker than the standing assumption of this paper:
\begin{equation}
\label{2.1.3}
f(u) \in C^2(\bar R_+) \,.
\end{equation}
Of course, radial solutions were studied prior to \cite{GNN}, in particular in a classical paper of D.D. Joseph and    T.S. Lundgren \cite{JL}, but 
the theorem of B. Gidas, W.-M. Ni and L. Nirenberg \cite{GNN} showed that radial solutions represent {\em all solutions} of the PDE (\ref{2.1.0}), which stimulated the interest in radial solutions.
\medskip

We now discuss the connection of the problem (\ref{2.1.0}) with Schr\"odinger's equation.
Let $v(x,t)$ be a complex-valued solution of a nonlinear Schr\"odinger's equation ($x \in R^n$, $t>0$)
\[
iv_t+\Delta v+v |v|^{p-1}=0 \,.
\]
Here $p>1$ is a constant, and $|v|$ denotes the complex modulus of $v$. Looking for the {\em standing waves}, one substitutes $v(x,t)=e^{imt}u(x)$, with a real valued $u(x)$, and a constant $m>0$. Then $u(x)$ satisfies
\[
\Delta u-mu+u |u|^{p-1}=0 \,.
\]
For a more general equation, where $v(x,t)$ is a complex-valued solution of
\[
iv_t+\Delta v+f(v)=0 
\]
a similar reduction works for any complex valued function $f(v)$, satisfying 
\[
f(e^{imt}u)=e^{imt} f(u) \,, \s \mbox{for any real $m$ and $u$} \,,
\]
and it leads to the equation
\[
\Delta u-mu+f(u)=0 \,.
\]

At first glance the equation (\ref{2.1.1}) appears to be  singular at $r=0$. In reality there is no singularity at $r=0$, after all this difficulty is introduced only by the spherical  coordinates. In fact, one can easily compute all of the derivatives $u^{(k)}(0)$, see P. Korman  \cite{K4}, then write down the Maclaurin series of solution, and show that it converges for small $r$, provided that $f(u)$ is analytic, see \cite{K4}. Since in this paper we only assume that $f(u) \in C^2(\bar R_+)$, let us recall the following result.

\begin{lma}\label{lma:2.1.1}
For any $\al >0$, consider  the initial value problem
\begin{eqnarray}
\label{2.1.5}
& u''(r)+\frac{n-1}{r} u'(r) +\la f(u(r))=0 \,, \\
& u(0)=\al, \, u'(0)=0 \,. \nonumber
\end{eqnarray}
Then one can find an $\epsilon >0$, so that this problem has a unique classical solution on the interval $[0,\epsilon)$.
\end{lma}
The proof can be found in  L.A. Peletier and J. Serrin \cite{PS1}, see also J.A. Iaia \cite{I}, or P. Quittner and P.  Souplet \cite{QS}.
\medskip 

We show next that the maximum value of solution $u(0)=\al$ can be used as a global parameter, i.e.,  it is impossible for two solutions of (\ref{2.1.1}) to share the same value of $\al$. This result has been known for a while, see e.g., E.N. Dancer \cite{D}.

\begin{lma}\label{lma:2.1.2}
The value of $u(0)=\al$ uniquely identifies the solution pair $(\la,u(r))$ of (\ref{2.1.1}) (i.e., there is at most one $\la$, with at most one solution $u(r)$, so that $ u(0)=\al$).
\end{lma}

\noindent {\bf Proof:} $\s$ 
Assume, on the contrary, that we have two solution pairs $(\la,u(r))$ and $(\mu,v(r))$, with $u(0)=v(0)=\al$. Clearly, $\la \ne \mu$, since otherwise we have a contradiction with uniqueness of initial value problems guaranteed by Lemma \ref{lma:2.1.1}.  (Recall that $u'(0)=v'(0)=0$.)
The change of variables $r=\frac{1}{\sqrt{\la}} t$ takes (\ref{2.1.5}) into
\begin{equation}
\label{2.1.5i}
 u''(t)+\frac{n-1}{t} u'(t)+f(u)=0, \s\s
u(0)=\al, \s u'(0)=0 \,. 
\end{equation}
The change of variables $r=\frac{1}{\sqrt{\mu}} t$ takes the equation for $v(r)$, at $\mu$, also into (\ref{2.1.5i}). By the preceding lemma, $u(t) \equiv v(t)$, but that is impossible, since $u(t)$ has its first root at $t=\sqrt{\la}$, while the first root of $v(t)$ is  $t=\sqrt{\mu}$.
\epf

We see that the value of $u(0)$ gives a {\em global parameter} on solution curves. In addition to its theoretical significance, this lemma is the key to numerical computation of bifurcation diagrams. For a sequence of values $\al=\al _i$ we compute the corresponding first root $r=r_i$ for  the solution of (\ref{2.1.5i}). (The initial value problem (\ref{2.1.5i}) is easy to solve numerically.) Then $\la _i=r_i^2$ is the corresponding value of $\la$ in (\ref{2.1.1}). We then plot all of the points $(\la _i,\al _i)$, obtaining the solution curve. According to the Lemma \ref{lma:2.1.2} above, these  two-dimensional curves give us a faithful picture of the solution set of the nonlinear PDE (\ref{2.1.1}).

\section{The shoot-and-scale method}
\setcounter{equation}{0}
\setcounter{thm}{0}
\setcounter{lma}{0}

\begin{figure}
\begin{center}
\scalebox{0.8}{\includegraphics{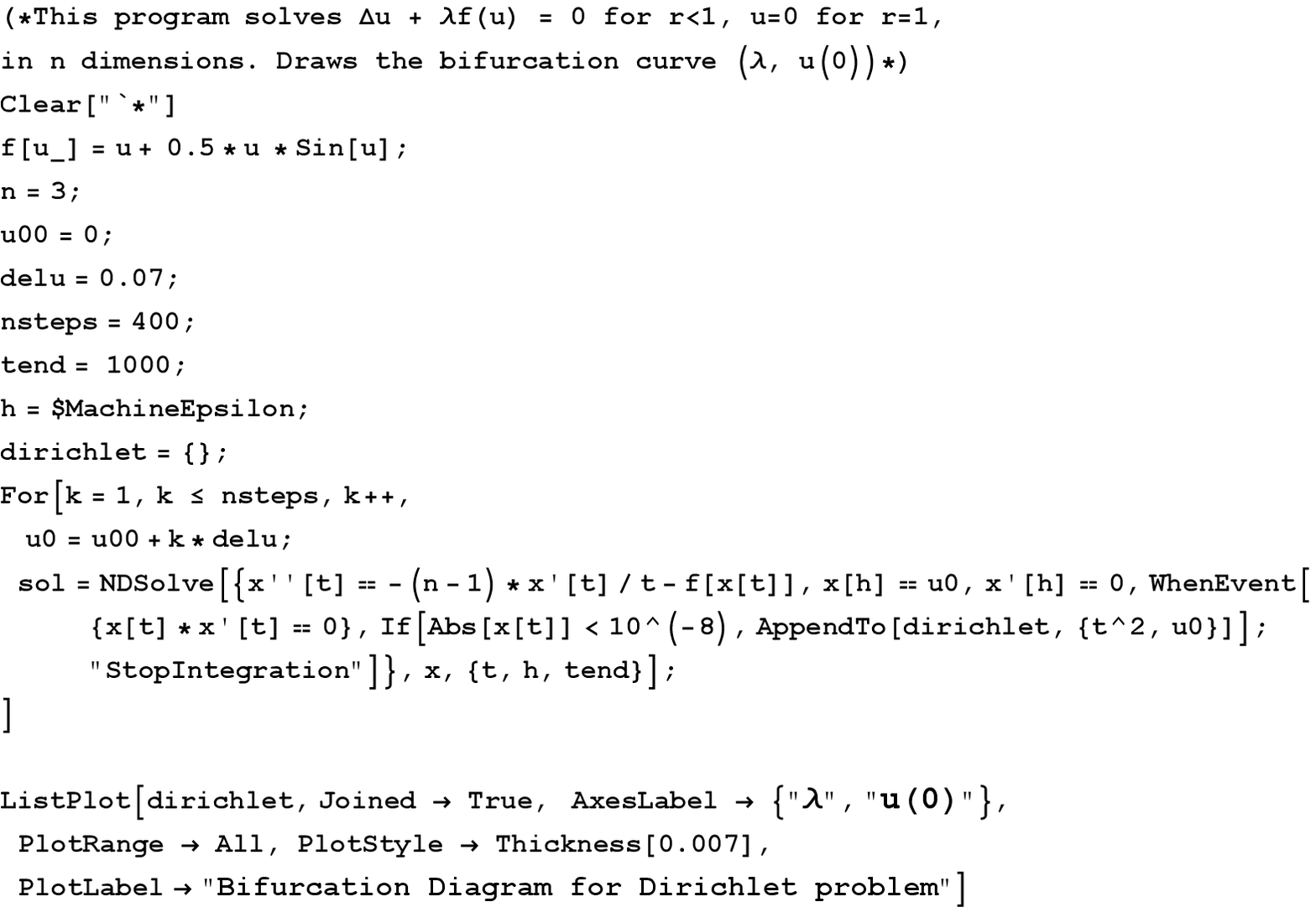}}
\end{center}
\caption{ The program to compute the  global solution curve for the problem (\ref{2.15.1}), with $f(u)=u+\frac{1}{2}u \sin u$, in space dimension $n=3$}
\label{web-c}
\end{figure}

Good analytical understanding of a problem goes hand in hand with efficient numerical calculation of its solutions. We know that for positive solutions the maximum value $u(0)=\al$ uniquely determines the solution pair $(\la,u(r))$ of the problem
\begin{equation}
\label{2.15.1}
\s\s u'' +\frac{n-1}{r}u'+\la f(u)=0, \s r \in (0,1), \s\s u'(0)=0, \; u(1)=0 \,,
\end{equation}
see Lemma \ref{lma:2.1.2} above. We also know that the parameter $\la $ in (\ref{2.15.1}) can be ``scaled out"
i.e., $v(r) \equiv u(\frac{1}{\sqrt {\la}} r)$ solves the equation 
\[
v''+\frac{n-1}{r}v'+f(v)=0 \,, 
\]
while $v(0)=u(0) \equiv \al$, and $v'(0)=u'(0)=0$. The  first root of $v(r)$ is $r=\sqrt {\la}$. We, therefore, solve the initial value problem
\begin{equation}
\label{2.15.2}
v''+\frac{n-1}{r}v'+f(v)=0, \s\s v(0)=\al, \s v'(0)=0,
\end{equation}
and compute its first positive root $r$. Then $\la =r^2$, by the above remarks. This way for each $\al$ we can find the corresponding $\la$.
After we choose sufficiently many $\al _n$'s and compute the corresponding $\la _n$'s, we  plot the pairs $(\la _n,\al _n)$, obtaining a bifurcation diagram in $(\la,\al)$ plane. If one needs to compute the actual solution $u(r)$ of (\ref{2.15.1}) at some point $(\la,u(0)=\al)$, this can be easily done by using the  NDSolve command in {\em Mathematica}. The program for producing a bifurcation diagram of  (\ref{2.15.1}), consists of  essentially one short loop.  It can be found on one of the authors web-page: http://homepages.uc.edu/$\sim$kormanp/, or  in the Figure \ref{web-c}.

\medskip

\begin{figure}
\begin{center}
\scalebox{0.6}{\includegraphics{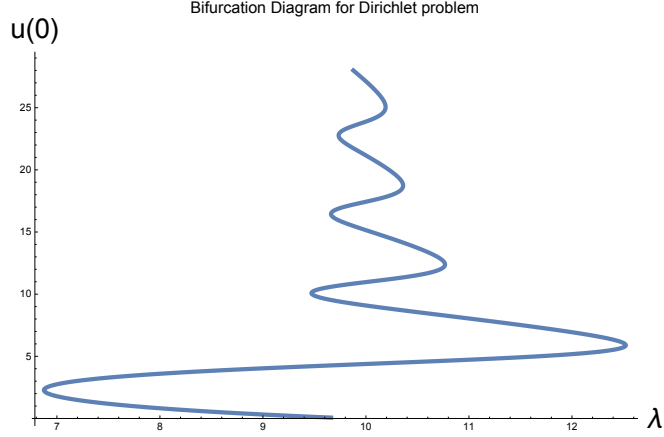}}
\end{center}
\caption{ The global solution curve for the problem (\ref{2.15.1}), with $f(u)=u+\frac{1}{2}u \sin u$, in space dimension $n=3$}
\label{web-c-pic}
\end{figure}

The shortness of this program may come as a surprise to some readers, considering that it produces the global solution curve for the problem (\ref{2.15.1}), with $f(u)=u+\frac{1}{2}u \sin u$. This bifurcation diagram is presented in Figure \ref{web-c-pic}. By the way, this computation illustrates the phenomenon of ``oscillatory bifurcation from infinity" first described  in   P. Korman \cite{K17}. Along the solution curve, $u(0) \ra \infty$, $\la \ra \la _1$ and $\la - \la _1$ changes sign infinitely many times, and moreover, $u(r)/u(0)$ tends to $ \varphi _1(r)$. Here $(\la _1, \varphi _1(r))$ is the principal eigenpair of the Dirichlet Laplacian on the unit ball, with $\varphi _1(0)=1$. What makes the situation different from the classical ``bifurcation from infinity" case is that $g(u)=\frac{1}{2}u \sin u$ term does not satisfy the condition $\lim _{u \ra \infty} \frac{g(u)}{u}=0$.
\medskip

We now describe the program. Beginning with $u(0) \equiv u00=0$, we choose the step size $delu=0.07$, and compute the $\la$'s corresponding to $u(0)=u00+n \, delu$, for the steps $n=1,2, \ldots, nsteps=400$. We implemented the shoot-and-scale method, using a highly sophisticated and accurate {\em Mathematica}'s command NDSolve to solve the initial value problem (\ref{2.15.2}). Formally, the equation has a singularity at $r=0$. We know that $u(0)$ is finite and given, and that $u'(0)=0$. Instead of starting the integration at zero and getting a zero by zero divide exception, we start at a value called  \$MachineEpsilon, which is a value near 0 and is  approximately $10^{-16}$.  We  integrate our equation (\ref{2.15.2}) until the first root of either $v(r)$ or $v'(r)$ is achieved. In case it is the first root of $v(r)$, the square of this root ($=\la$) and the corresponding value of $v(0)=u(0)$  are stored in the file named ``dirichlet", which is  later plotted. {\em Mathematica} requires us to specify the upper limit of integration, so we choose a large number $tend=1000$, which is never achieved in practice.
\medskip

A simple adjustment is needed if the solution curve has  disjoint  branches. Replace the option Joined $\ra$ True by  Joined $\ra$ False in the ListPlot command, to avoid the branches being joined. A nicer picture can be obtained by first plotting each branch separately, and then plotting them jointly (using    {\em Mathematica}'s command Show).

\section{The shoot-and-scale method in case of $p$-Laplacian}
\setcounter{equation}{0}
\setcounter{thm}{0}
\setcounter{lma}{0}

We now apply the shoot-and-scale method to compute curves of positive solutions for the Dirichlet problem  in case of the $p$-Laplacian
\begin{equation}
\label{sas1}
\p (u')'+\frac{n-1}{r} \, \p (u') +\lambda f(u)=0, \s  u'(0)=u(1)=0 \,,
\end{equation}
where $\p (v)=v |v|^{p-2}$, $p>1$. As in the Laplacian case $p=2$, the value of $u(0)$ uniquely identifies the solution pair $(\lambda,u(r))$, so that the solution curves can be drawn in  the $(\lambda, u(0))$ plane. The proof is exactly the same as in case $p=2$, using the following existence and uniqueness result from P. Korman \cite{K4}.
\begin{lma}
Assume that $f(u)$ is Lipschitz continuous. Then one can find an $\epsilon >0$ so that the problem
\[
\p (u')'+\frac{n-1}{r} \, \p (u') +\lambda f(u)=0 \,, \s u(0)=\alpha \,, \; u'(0)=0
\]
has a unique solution $u(r) \in C^1[0,\epsilon ) \cap C^2(0,\epsilon )$, for any $\alpha >0$.
\end{lma}
\medskip

We have $\p (u')'=\p' (u') u''=(p-1) |u'|^{p-2}u''$, using that $\p' (v)= \newline (p-1) |v|^{p-2}$. Dividing (\ref{sas1}) by $(p-1)|u'|^{p-2}$,  we have 
\begin{equation}
\label{sas2}
u''+\frac{n-1}{(p-1)r}\, u'+\lambda \, \frac{f(u)}{(p-1) |u'|^{p-2}}=0, \s  u'(0)=u(1)=0 \,.
\end{equation}
Similarly to the case of the Laplacian ($p=2$), one may  approach (\ref{sas2}) by solving the initial value problem
\begin{equation}
\label{sas3}
u''+\frac{n-1}{(p-1)r}\, u'+ \frac{f(u)}{(p-1) |u'|^{p-2}}=0, \s  u(0)=\al >0, \s u'(0)=0 
\end{equation}
until its solution $u(r)$ reaches the first root at some $\xi$. (The form (\ref{sas3}) of our equation is convenient for {\em Mathematica}'s NDSolve command to solve.) Then by scaling we obtain a solution of (\ref{sas2}) at $\lambda =\xi ^p$, and a point $(\lambda,\al )$ on the solution curve in the $(\lambda, u(0))$ plane. After computing a sufficient number of such points, one  plots the solution curve. Observe that for $p>2$ we have a {\em true singularity at $r=0$} in the last term of (\ref{sas3}) (in case $p=2$ there is only an apparent singularity at $r=0$). It was proved in P. Korman \cite{K2015} that the solution of  (\ref{sas3}) satisfies 
\beq
\lbl{sas3a}
u(r) \approx \alpha +a_1 r^{\frac{p}{p-1}} \,, \s \mbox{ for $r$ small} \,,
\eeq
where the constant $a_1<0$ is computed using the formula (3.3) of that paper. A simple-minded approach involves choosing a small $h$, say $h=0.00001$, then approximating $u(h)$ and  $u'(h)$ using (\ref{sas3a}), and use these values as the initial conditions while integrating (\ref{sas3}) for $r>h$. Remarkably, this approach works reasonably well, even for fast growing $f(u)$ like  $f(u)=e^u$, even though the last term in (\ref{sas3}) is huge at $r=h$. This is due to the high accuracy of the {\em Mathematica}'s NDSolve command. 
\medskip

A better approach is to use the {\em regularizing transformation}  $r \ra z$, given by 
\beq
\lbl{sas5}
z=r^{\frac{p}{2(p-1)}} \,,
\eeq
which was introduced in P. Korman \cite{K2015}. We present next a simplified proof of that result.

\begin{lma}\lbl{lma:1}
Denote 
\beq
\lbl{22a}
\bar \beta =\frac{p}{2(p-1)}, \s a=\bar \beta ^p \,, \s A=\beta ^{p-1} (\beta -1)+\frac{n-1}{p-1}  \beta ^{p-1} \,.
\eeq
Then, for $p>2$, the change of variables (\ref{sas5}) transforms (\ref{sas3}) into
\beqa
\lbl{333}
& au''(z)+\frac{A}{z}u'(z)+\frac{z^{p-2}}{(p-1)|u'(z)|^{p-2}} f(u)=0 \,, \\ 
& u(0)=\al ,   \; u'(0)=0 \,. \nonumber
\eeqa 
\end{lma}

\pf
Consider the transformation $z=r^{\beta} $ in  (\ref{sas3}), with the constant $\beta$ to be specified. Using that $u_r=\beta r^{\beta -1}u_z$ and $u_{rr}=\beta ^2 r^{2\beta -2}u_{zz}+\beta (\beta-1)r^{\beta -2} u_z$ obtain from (\ref{sas3})
\[
\beta ^2 r^{2\beta -2}u_{zz}+\left(\beta (\beta-1)+\beta \frac{n-1}{p-1} \right) r^{\beta -2} u_z+\frac{f(u)}{(p-1) \beta ^{p-2} r^{(\beta -1)(p-2)} |u_z|^{p-2}}=0 \,.
\]
Divide this equation by $r^{2\beta -2}$,  and multiply by $\beta ^{p-2}$ to obtain
\[
au''(z)+\frac{A}{z}u'(z)+\frac{1}{(p-1) r^{(\beta -1)p}|u'(z)|^{p-2}} f(u)=0 \,.
\]
This equation becomes (\ref{333}), provided we choose $\beta$ so that
\[
\frac{1}{r^{(\beta -1)p}}=z^{p-2}=r^{\beta (p-2)} \,,
\]
which occurs for $\beta=\bar \beta =\frac{p}{2(p-1)}$.
\epf

The equation in (\ref{333}) has no singularity at $z=0$, in fact the solution is of class $C^2$ for continuous $f(u)$ (all of the derivatives $u^{(n)}(0)$ can be computed if $f(u) \in C^{\infty}$). It follows that solutions of (\ref{sas1}) are of the form  $u(r)=w \left(r^{\frac{p}{2(p-1)}} \right)$, where the function $w(z)$ is twice differentiable at $z=0$, $w(z) \in C^2[0,1)$. However, for the {\em Mathematica}'s NDSolve command the equation (\ref{333})  is still considered to be singular. Therefore we use the above mentioned approximation from P. Korman \cite{K2015}
\beq
\lbl{sas6}
u(z) \approx \al+a_1 z^2+a_2 z^4 \,, 
\eeq
to approximate $u(h)$ and $u'(h)$, and then perform the numerical calculation for $z>h$. Here $a_1<0$ is given by the formula (3.3) of \cite{K2015}:
\[
a_1=-\left[ \frac{1}{(p-1)2^{p-2}B_1} f(\al ) \right]^{\frac{1}{p-1}} \,,
\]
where $B_1=2a+\frac{2A}{p-1}$.  We now calculate $a_2$ using the formula (3.4) of \cite{K2015}:
\[
a_2=-\frac{1}{B_2 C_2} \lim _{z \ra 0} \frac{\frac{z^{p-2}}{(p-1) \left(2 a_1 z \right)^{p-2}} f \left( \al +a_1 z^2 \right)+2a_1a+\frac{2a_1A}{p-1}}{z^2} \,.
\]
This limit does not change if we replace $f \left( \al +a_1 z^2 \right)$ by $f(\al)+f'(\al) a_1z^2$, which leads to
\[
a_2=-\frac{f'(\al) a_1}{B_2 C_2(p-1) \left(2a_1 \right)^{p-2}}  \,.
\]
If $z_0$ denotes the first root of $u(z)$, while $\xi$ is the first root of $u(r)$, then in view of (\ref{sas5}) the value of $\la$ corresponding to $u(0)=\al$ is 
\[
\la =\xi ^p=z_0^{2(p-1)} \,.
\]

\begin{figure}
\begin{center}
\scalebox{0.65}{\includegraphics{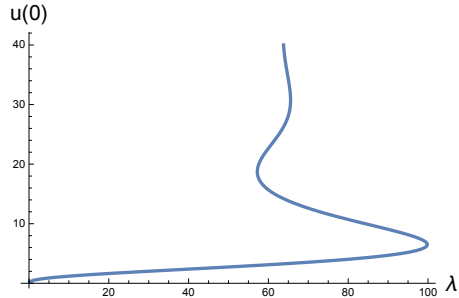}}
\end{center}
\caption{ The   global  solution curve for the $p$-Laplace problem (\ref{sas1}), with $f(u)=e^u$, $p=4$ and $n=5$}
\label{plap}
\end{figure}

After calculating sufficiently many points $(\la _n, \al _n)$, we plot the solution curve for the problem (\ref{sas1}). An example is presented in Figure \ref{plap}.

\section{Bifurcation diagrams for Neumann boundary conditions}
\setcounter{equation}{0}
\setcounter{thm}{0}
\setcounter{lma}{0}

A small modification of the above program produces positive decreasing solutions of the Neumann 
problem
\beqa
\label{10}
& u'' +\frac{n-1}{r}u'+\la f(u)=0 \,, \s r \in (0,1), \s\s u'(0)= u'(1)=0 \\ \nonumber
& u(r)>0 \,, \s u'(r)<0 \,, \s r \in (0,1) \,. \nonumber
\eeqa
There is no analog of B. Gidas, W.-M. Ni and L. Nirenberg \cite{GNN} result for Neumann 
problem, but radial solutions on a ball are still of interest. 
\medskip

\begin{figure}
\begin{center}
\scalebox{0.8}{\includegraphics{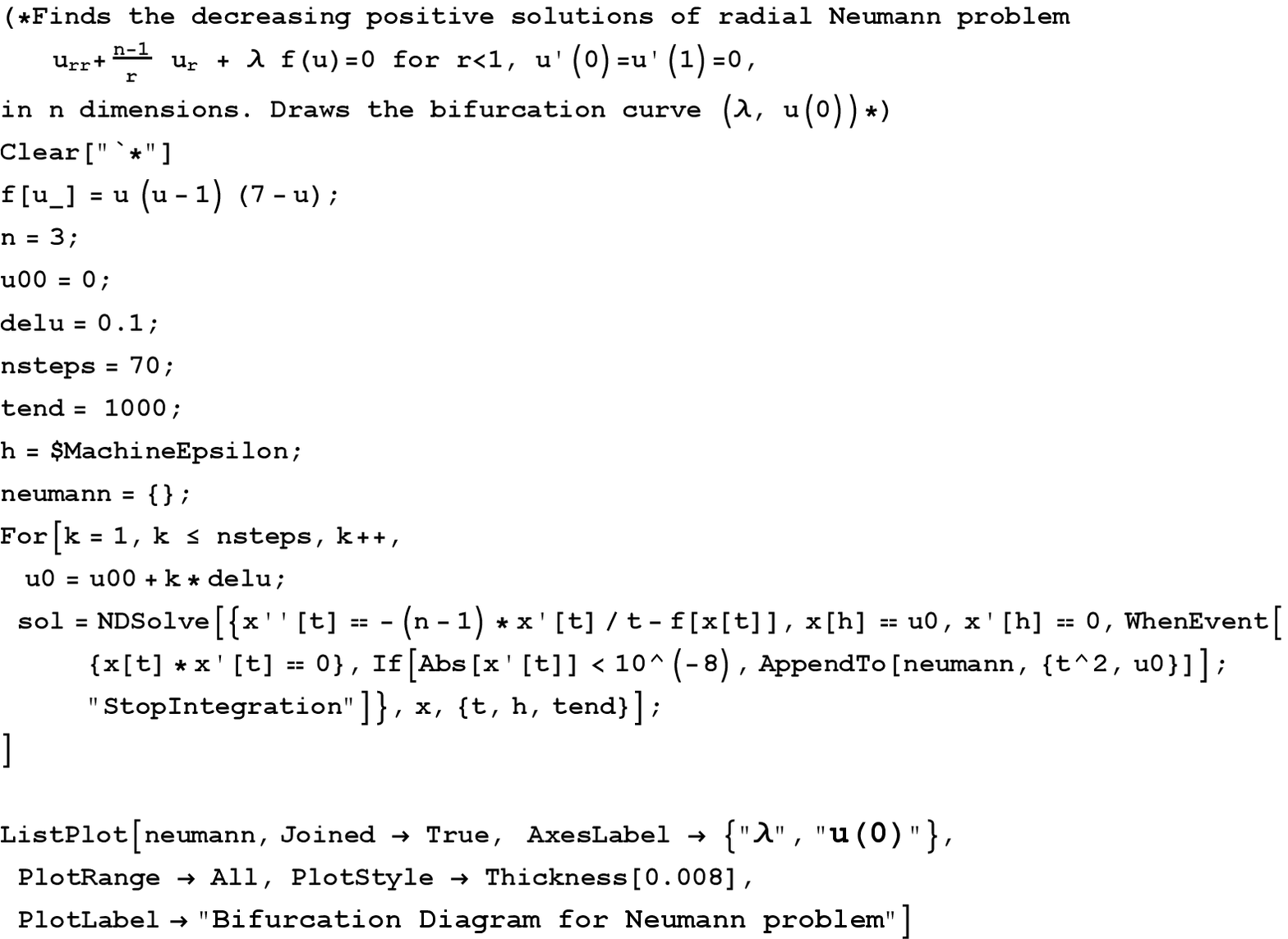}}
\end{center}
\caption{ The program to compute the  global Neumann solution curve for the problem (\ref{10}), with $f(u)=u(u-1)(7-u)$}
\label{neumann}
\end{figure}

\begin{figure}
\begin{center}
\scalebox{0.6}{\includegraphics{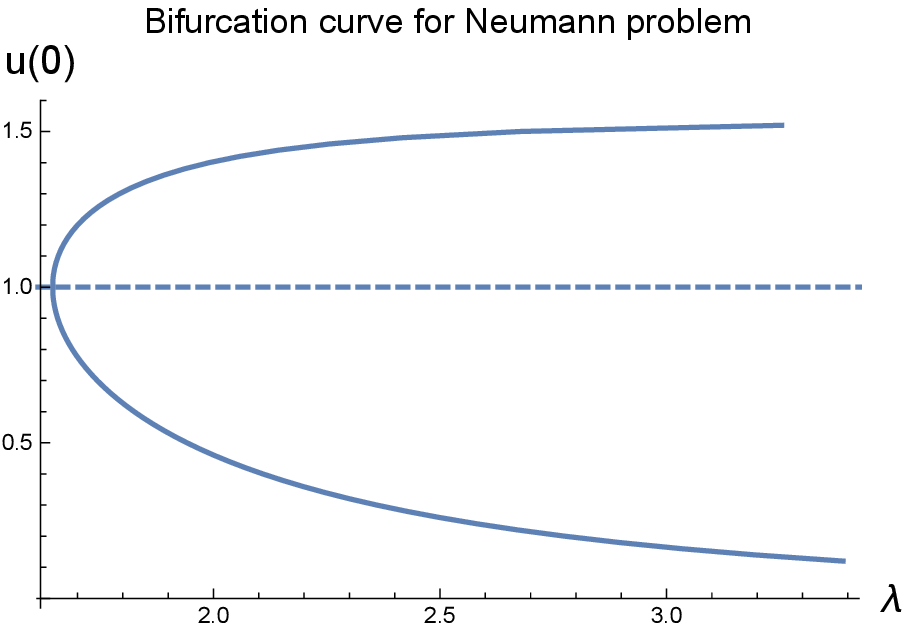}}
\end{center}
\caption{ The   global Neumann solution curves for the problem (\ref{10}), with $f(u)=u(u-1)(7-u)$, $n=1$}
\label{neumann-picture-1}
\end{figure}

\begin{figure}
\begin{center}
\scalebox{0.6}{\includegraphics{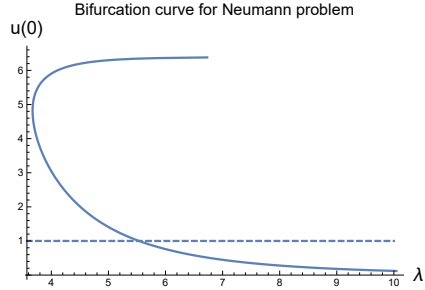}}
\end{center}
\caption{ The global Neumann solution curves for the problem (\ref{10}), with $f(u)=u(u-1)(7-u)$, $n=5$}
\label{neumann-picture-2}
\end{figure}

In Figure \ref{neumann} we present the program to compute the bifurcation diagrams for (\ref{10}), in case $f(u)=u(u-1)(7-u)$ (electronically available at http://homepages.uc.edu/$\sim$kormanp/). For this $f(u)$ the problem (\ref{10}) has a trivial solution $u(r)=1$ for all $\la$. Our computations show that the situation is different for the cases $n=1$ and $n \geq 2$. In the case $n=1$ there are two half-branches bifurcating from $u=1$ at some $\la _0 \approx 1.65$ that continue for all $\la > \la _0$ without any turns, see Figure \ref{neumann-picture-1}. That both half-branches bifurcating from $u=1$ do not turn was proved rigorously by R. Schaaf \cite{S10}, in case $n=1$, for a class of nonlinearities $f(u)$ which includes polynomials with simple roots. In Figure \ref{neumann-picture-2} we present the bifurcation picture for $n=5$, which is typical for all other dimensions $n \geq 2$ that we tried. This time two half-branches bifurcate from $u=1$ at some $\la _1 \approx 5.5$. The lower branch continues without any turns for all $\la >\la _1$, while the upper branch continues for decreasing $\la$, and makes exactly one turn. Proving these facts rigorously is a  very challenging open problem. The bifurcation methods as in e.g., P. Korman \cite{K2} or T. Ouyang and J. Shi \cite{OS} appear to be inapplicable here. 

\medskip

We now describe the program presented in  Figure \ref{neumann}. We integrate the initial value problem (\ref{2.15.2}), similarly to the Dirichlet case. Integrations are stopped if either the root of $v(r)$ or of $v'(r)$ is  reached. We record these values of $r$ and  the corresponding $v(0)$ in two separate files, and then plot the file involving $v'(r)=0$, obtaining the curve of positive decreasing solutions of the Neumann problem (\ref{10}).

\section{Solution curves for non-autonomous problems }
\setcounter{equation}{0}
\setcounter{thm}{0}
\setcounter{lma}{0}

We wish to compute the curves of  positive solutions for non-autonomous problems 
\begin{equation}
\label{4a.9}
\Delta u +\la f(|x|,u)=0 \;\; \mbox{in $B$}, \;\; u=0 \s \mbox{on $\partial B$} \,,
\end{equation}
where $B$ is the unit ball $|x|<1$ in $R^n$. We shall assume that the function $f(r,u) \in C^2(B \times R_+)$ satisfies
\begin{equation}
\label{4a.10}
f(r,u)>0,  \s \s \mbox{for $0<r<1$ and $u>0$} \,,
\end{equation}
\begin{equation}
\label{4a10.1}
 f_r(r,u) \leq 0, \s \s \mbox{for $0<r<1$ and $u>0$} \,.
\end{equation}
By B. Gidas, W.-M. Ni and L. Nirenberg \cite{GNN}, any positive solution of (\ref{4a.9}) is radially symmetric, i.e., $u=u(r)$, $r=|x|$, and hence it satisfies
\begin{equation}
\label{4a.11}
u''+\frac{n-1}{r}u'+\la f(r,u)=0, \s \s  u'(0)=u(1)=0 \,.
\end{equation}
The shoot-and-scale method does not apply for non-autonomous problems (one can still shoot, but the scaling fails.)
\smallskip

The standard approach to numerical computation of solutions involves curve following, i.e., continuation in $\la$ by using the predictor-corrector type methods, see e.g., E.L. Allgower and K. Georg \cite{A}. These methods are well developed, but not easy to implement, as the solution curve $u=u(x,\la)$ may consist of several parts, each having multiple turns. Here $\la $ is a local parameter, but usually it is not a global one, because of the possible turning points. We wish to perform the continuation using a global parameter, similarly to the shoot-and-scale method.
\smallskip

The quantity $\al = u(0)$ gives the maximum value of any positive solution. In case $n=1$,  the value of $\al = u(0)$ is known to be  a global parameter, i.e., the value of $\al$ uniquely identifies the solution pair $(\la,u(r))$, see P. Korman \cite{K2013} or P. Korman and J. Shi \cite{KS1}. In case $n>1$, we shall still do continuation in $\al$,  computing the solution curve of (\ref{4a.11}) in the form $\la=\la (\al)$, although now there is no guarantee that all positive solutions of (\ref{4a.11}) are obtained. We begin with a simple lemma.

\begin{lma}\label{lma:n1}
The solution of the linear problem
\[
u'' +\frac{n-1}{r}u'+  g(r)=0 \s\s \mb{for $0<r<1$,} \s\s u(0)=\al, \s u'(0) =0 
\]
can be represented in the form
\[
u(r)=\al+\frac{1}{n-2}r^{-n+2} \int_0^r \left(z^{n-2}-r^{n-2} \right) \, z g(z) \, dz, \s\s \mbox{for $n \ne 2$} \,,
\]
\[
u(r)=\al+ \int_0^r \left(\ln z -\ln r \right) \, z g(z) \, dz, \s\s \mbox{for $n = 2$} \,.
\]  
\end{lma}

\noindent {\bf Proof:} $\s$
Integrating
\[
\left(r^{n-1}u'(r) \right)'=-r^{n-1} g(r)
\]
over the interval $(0,z)$, we express
\[
u'(z)=-\frac{1}{z^{n-1}} \int _0^z t^{n-1} g(t) \, dt \,.
\]
Integrating over the interval $(0,r)$, we have
\[
u(r)=\al- \int_0^r \frac{1}{z^{n-1}} \int _0^z t^{n-1} g(t) \, dt \,dz\,.
\]
Integrating by parts in the last integral (with $u=\int _0^z t^{n-1} g(t) \, dt$, $dv=\frac{1}{z^{n-1}} \,dz$ ), we conclude the proof.
\hfill$\diamondsuit$ \medskip

If we solve the initial value problem
\begin{equation}
\label{4a.n3}
u'' +\frac{n-1}{r}u'+ \la f(r,u)=0, \s\s u(0)=\al,  \s\s u'(0)=0  \,,
\end{equation}
then we need to find $\la$, so that $u(1)=0$, in order to obtain a solution of (\ref{4a.11}). By Lemma \ref{lma:n1}, we rewrite the equation (\ref{4a.11}) in the integral form (for $n  \ne 2$)
\[
u(r)=\al+\frac{\la}{n-2}r^{-n+2} \int_0^r \left(z^{n-2}-r^{n-2} \right) \, z f(z,u(z)) \, dz, \s\s \mbox{for $n \ne 2$} \,,
\]
and then the equation for $\la$ is
\begin{equation}
\label{4a.n4}
F(\la) \equiv u(1)= \al+\frac{\la}{n-2} \int_0^1 \left(z^{n-2}-1 \right) \, z f(z,u(z)) \, dz =0 \,.
\end{equation}
We solve this equation by using Newton's method
\[
\la _{n+1}=\la _{n}-\frac{F(\la _{n})}{F \,'(\la _{n})} \,, \s n=0,1,2, \ldots \,.
\]
The iterations begin at $\la _0$, which we choose to be the value of $\la$ corresponding to the preceding value of $\al =u(0)$.
Calculate
\[
F(\la _{n})=\al +\frac{\la_n}{n-2} \int_0^1 \left(z^{n-2}-1 \right) \, z f(z,u(z)) \, dz \,,
\]
\[
F'(\la _{n})=\frac{1}{n-2} \int_0^1 \left(z^{n-2}-1 \right) \, z f(z,u(z)) \, dz   
\]
\[
+\frac{\la_n}{n-2} \int_0^1 \left(z^{n-2}-1 \right) \, z f_u(z,u(z)) u_{\la} \, dz \,,
\]
where $u=u(r,\la _n)$ and $u_{\la}=u_{\la}(r,\la _n)$ are respectively the solutions of 
\begin{equation}
\label{4a.n5}
u'' +\frac{n-1}{r}u' + \la _n f(r,u)=0, \s u(0)=\al,  \s u'(0)=0  \,,
\end{equation}
and
\begin{equation}
\label{4a.n6}
\s\s\s  u_{\la}''+\frac{n-1}{r}u_{\la}' + \la _n  f_u(r,u(r,\la _n))u_{\la}+f(r,u(r,\la _n))=0
\end{equation}
\[
 \s u_{\la}(0)=0,  \; u_{\la}'(0)=0  \,.
\]
(As we vary $\la$, when solving (\ref{4a.n4}), we keep $u(0)=\al$ fixed, and that is the reason why $u_{\la}(0)=0$.)
This method is very easy to implement. It  requires  repeated solutions of the initial value problems (\ref{4a.n5}) and (\ref{4a.n6}) (using the NDSolve command in {\em Mathematica}).
\medskip

In case $n=2$, we have
\[
F(\la)=\al + \la \int_0^1  z \ln z f(z,u(z)) \, dz \,,
\]
\[
F'(\la )=\int_0^1  z \ln z f(z,u(z)) \, dz+\la \int_0^1  z \ln z f_u(z,u(z)) u_{\la}(z) \, dz \,,
\]
and the rest is as before.
\medskip

\begin{figure}
\begin{center}
\scalebox{0.6}{\includegraphics{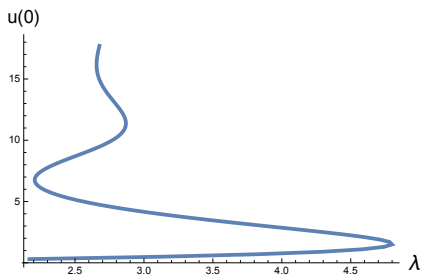}}
\end{center}
\caption{ The  global solution curve for  the  problem (\ref{4a.n7}), $n=3$}
\label{gelfand-gr1}
\end{figure}

\begin{figure}
\begin{center}
\scalebox{0.6}{\includegraphics{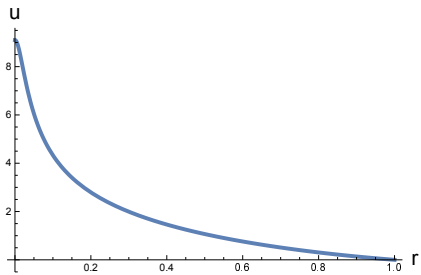}}
\end{center}
\caption{ The  solution $u(r)$, corresponding to $\al=u(0)=9.1$, and $\la \approx 2.59566$}
\label{gelfand-gr2}
\end{figure}

\noindent
{\bf Example} $\s\s$ We  solved the problem 
\begin{eqnarray}
\label{4a.n7}
& u''+\frac{n-1}{r}u' + \la \left(1-1.1r^2 \right)e^u=0 \,, \s\s \mb{for $0<r<1$,} \\ \nonumber
& u'(0)=u(1)=0  \,, \nonumber
\end{eqnarray}
for the Gelfand's equation of combustion theory, with a sign-changing potential $1-1.1r^2$.
The global curve of positive solutions, for $n=3$, is presented in Figure \ref{gelfand-gr1}. For any point $(\la, \al)$ on this curve ($\al =u(0)$), the actual solution $u(r)$ is easily computed by shooting (using the NDSolve command in {\em Mathematica}), i.e., by solving (\ref{4a.n3}). In Figure \ref{gelfand-gr2} we present the solution $u(r)$ for $\la \approx 2.59566$, when $u(0)=9.1$. (This solution lies on our solution curve in Figure \ref{gelfand-gr1}, after the second turn.)
The picture in Figure \ref{gelfand-gr1} suggests that solution curve for the problem (\ref{4a.n7}) has infinitely many turns, similarly to the case of the constant potential, according to the classical result of D.D. Joseph and    T.S. Lundgren \cite{JL}.

\section{Solution curves for the elastic beam equation}
\setcounter{equation}{0}
\setcounter{thm}{0}
\setcounter{lma}{0}

We now compute  curves of positive solutions  for  the  equation modeling the displacements of an elastic beam, which is clamped at both end points,
\begin{eqnarray} 
\lbl{bm1}
& u''''(x)=\la f(u(x)), \s  \mbox{for $x \in (-1,1)$} \\ \nonumber
&  u(-1)=u'(-1)=u(1)=u'(1)=0 \,, \nonumber
\end{eqnarray}
obtaining the solution curves in $(\lambda,u(0))$ plane, similarly to the second order equations. Under the assumption that
\beq
\lbl{bm2}
f(u)>0 \,, \s \mbox{and} \s f'(u)>0 \s\s \mbox{for $u>0$}
\eeq
it was shown in P. Korman \cite{K} that any positive solution of (\ref{bm1}) is an even function taking its global maximum at $x=0$. Moreover, the value of $u(0)$ is a global parameter, uniquely identifying the solution pair $(\la,u(x))$  (see also P. Korman and J. Shi \cite{KS} for a generalization). That $u(0)$ is a global parameter is quite a surprising result, because unlike the second order equations, the knowledge of $u(0)=\al$ and $u'(0)=0$ does not identify the solutions of the initial value problem corresponding to the equation in (\ref{bm1}) (we have $u'''(0)=0$ by the symmetry of solutions, but one also needs to know the value of $u''(0)=\beta$, in order to ``shoot").
\medskip

As we just mentioned, the initial conditions corresponding to the positive symmetric solutions of (\ref{bm1}) are
\beq
\lbl{bm3}
u(0)=\al >0 \,, \s u'(0)=0 \,, \s u''(0) =\beta \leq 0 \,, \s u'''(0)=0 \,.
\eeq
We record the following simple observation.

\begin{lma}
For a given continuous function $g(x)$ the solution of 
\[
u''''=g(x) \,,
\]
together with the initial conditions (\ref{bm3}) is 
\beq
\lbl{bm4}
u(x)=\int_0^x \frac{(x-t)^3}{3!} \, g(t) \, dt+\al+\frac{\beta x^2}{2} \,.
\eeq
\end{lma}

We know that the value of $u(0)=\al$ uniquely identifies the solution pair $(\la,u(x))$ (and plotting of $\la=\la(\al )$ gives the bifurcation diagram). However, we do not know the corresponding value of $u''(0)=\beta$, which makes it impossible to shoot and scale. We now describe an algorithm to compute both $(\beta,\la)$, given $\al$.
\medskip

We begin by converting the equation in (\ref{bm1}), together with the initial conditions (\ref{bm3}) into an integral equation
\[
u(x)=\la \int_0^x \frac{(x-t)^3}{3!} \, f(u(t)) \, dt+\al+\frac{\beta x^2}{2} \,,
\]
by using (\ref{bm4}). To obtain the solution of our problem (\ref{bm1}) corresponding to $u(0)=\al$, we need to find $(\beta,\la)$ so that
\beq
\lbl{bm0}
u(1)=\la \int_0^1 \frac{(1-t)^3}{3!} \, f(u(t)) \, dt+\al+\frac{\beta }{2}=0 \,,
\eeq
\[
u'(1)=\la \int_0^1 \frac{(1-t)^2}{2} \, f(u(t)) \, dt+\beta=0 \,.
\]
Denoting 
\[
F(\la ,\beta)=\la \int_0^1 \frac{(1-t)^3}{3!} \, f(u(t)) \, dt+\al+\frac{\beta }{2} \,, \s G(\la,\beta)=\la \int_0^1 \frac{(1-t)^2}{2} \, f(u(t)) \, dt+\beta \,, 
\]
we recast the  equations (\ref{bm0}) as
\beqa
\lbl{bm5}
& F(\la ,\beta)=0 \\ \nonumber
& G(\la,\beta)=0 \,.\nonumber
\eeqa
We use Newton's method to solve the system (\ref{bm5}). Assuming that the $n$-th iterate $(\la _n,\beta _n)$ is already computed, we linearize (\ref{bm5}) at this point
\beqa
\lbl{bm6}
& F_{\la}(\la _n,\beta _n)(\la -\la _n)+F_{\beta}(\la _n,\beta _n)(\beta -\beta _n)=0 \\ \nonumber
& G_{\la}(\la _n,\beta _n)(\la -\la _n)+G_{\beta}(\la _n,\beta _n)(\beta -\beta _n)=0 \,,\nonumber
\eeqa 
then solve this system for $(\la ,\beta)$, and declare the solution to be our next iterate $(\la _{n+1},\beta _ {n+1})$. Observe that the solution $u(t)$ of (\ref{bm1}) and (\ref{bm3}), which was used in the definitions of $F(\la ,\beta)$ and $G(\la ,\beta)$ is not known precisely, so we shall replace it by the best approximation available at each iteration step. Namely, define $u_n(x)$ to be the solution of the following initial value problem
\beqa
\lbl{bm7}
& u''''(x)=\la _n f(u(x)) \\ \nonumber
& u(0)=\al \,, \s u'(0)=0 \,, \s u''(0) =\beta _n  \,, \s u'''(0)=0 \,. \nonumber
\eeqa
Clearly, the solution $u(x)$ of (\ref{bm1}) and (\ref{bm3}) depends on both $\la$ and $\beta$. We denote by $u_{\la}(x)$ the derivative of $u(x)$ with respect to $\la$, and by $u_{\beta}(x)$ the derivative of $u(x)$ with respect to $\beta$. Then we compute the quantities in (\ref{bm6}) as follows
\[
F_{\la}(\la _n,\beta _n)= \int_0^1 \frac{(1-t)^3}{3!} \, f(u_n(t)) \, dt+\la _n \int_0^1 \frac{(1-t)^3}{3!} \, f'(u_n(t)) u_{\la} (t) \, dt \,,
\]
\[
F_{\beta}(\la _n,\beta _n)=\la _n \int_0^1 \frac{(1-t)^3}{3!} \, f'(u_n(t)) u_{\beta} (t) \, dt +\frac{1}{2} \,,
\]
\[
G_{\la}(\la _n,\beta _n)= \int_0^1 \frac{(1-t)^2}{2} \, f(u_n(t)) \, dt+\la _n \int_0^1 \frac{(1-t)^2}{2} \, f'(u_n(t)) u_{\la} (t) \, dt \,,
\]
\[
G_{\beta}(\la _n,\beta _n)=\la _n \int_0^1 \frac{(1-t)^2}{2} \, f'(u_n(t)) u_{\beta} (t) \, dt +1 \,.
\]
To compute $u_{\la} (x)$ one needs to solve the following linear initial value problem
\beqa
\lbl{bm8}
& u_{\la}''''(x)= f(u_n(x))+\la _n f'(u_n(x)) u_{\la}(x)\\ \nonumber
& u_{\la}(0)=0 \,, \s u_{\la}'(0)=0 \,, \s u_{\la}''(0) =0  \,, \s u_{\la}'''(0)=0 \,. \nonumber
\eeqa
obtained by differentiating (\ref{bm1}) and (\ref{bm3}) in $\la$. (We have $u_{\la}(0)=0$, because $u(0)$ is kept fixed at $\al$, and $u_{\la}''(0) =0$, since only $\la$ is varied, and $\beta$ is kept fixed.) To compute $u_{\beta} (x)$ one needs to solve the following linear initial value problem
\beqa
\lbl{bm9}
& u_{\beta}''''(x)= \la _n f'(u_n(x)) u_{\beta}(x)\\ \nonumber
& u_{\beta}(0)=0 \,, \s u_{\beta}'(0)=0 \,, \s u_{\beta}''(0) =1  \,, \s u_{\beta}'''(0)=0 \,. \nonumber
\eeqa
obtained by differentiation of (\ref{bm1}) and (\ref{bm3}) in $\beta$. 
\medskip

We now review the computation of the global curve $\la =\la ( \al)$ of positive solutions for the problem (\ref{bm1}). We choose $\al _0 $, the initial value for $\al$, the step-size called $step$, and the number of steps $N$, then for $\al _m=\al _0+m \, step$, $m=1,2, \ldots N$, we calculate the corresponding values of $\la$ by obtaining $(\la, \beta)$ as a limit of Newton's iterates $(\la _n,\beta _n)$. Given  $(\la _n,\beta _n)$, the next Newton's iterate $(\la _{n+1},\beta _{n+1})$ is produced as follows.

\begin{enumerate}
  \item Solve (\ref{bm7}) to calculate $u_n(x)$\,,
  \item Solve (\ref{bm8}) to calculate $u_{\la}(x)$\,,
  \item Solve (\ref{bm9}) to calculate $u_{\beta}(x)$\,,
  \item Calculate the numbers $F_{\la}(\la _n,\beta _n)$, $F_{\beta}(\la _n,\beta _n)$, $G_{\la}(\la _n,\beta _n)$ and $G_{\beta}(\la _n,\beta _n)$\,,
  \item Solve the linear system (\ref{bm6}) for $(\la,\beta)$ which is the desired $(\la _{n+1},\beta _{n+1})$.
  \end{enumerate}
Newton's method is started at $(\la _0, \beta _0)$, for which we use the values of $\la $ and $\beta$ at the preceding value of $\al$, namely $\al=\al _{m-1}$.
\medskip

\begin{figure}
\begin{center}
\scalebox{0.65}{\includegraphics{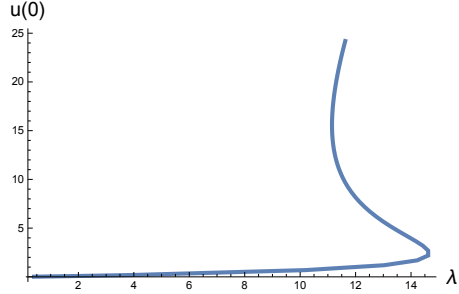}}
\end{center}
\caption{ The  global solution curve for  the  problem (\ref{bm20})}
\label{fourth}
\end{figure}

\noindent
{\bf Example} $\s\s$ We  used {\em Mathematica} to compute the solutions of  the problem 
\begin{eqnarray} 
\lbl{bm20}
& u''''(x)=\la e^{\frac{au}{a+u}}, \s  \mbox{for $x \in (-1,1)$} \\ \nonumber
&  u(-1)=u'(-1)=u(1)=u'(1)=0 \,, \nonumber
\end{eqnarray}
with $a=5$.
The $S$-shaped global curve of positive solutions is presented in Figure \ref{fourth}. For any point $(\la, \al)$ on this curve, the actual solution $u(x)$ is easily computed by shooting (using the NDSolve command in {\em Mathematica}), i.e., by solving the equation in (\ref{bm1}) together with the initial conditions (\ref{bm3}), using the function $\beta=\beta(\al )$ which was calculated while producing the above bifurcation diagram. For second order equations, the function $f(u)= e^{\frac{au}{a+u}}$ is prominent in connection to the ``perturbed Gelfand equation" of combustion theory, and the solution curve is also $S$-shaped, see e.g., \cite{K2} for the references.

\medskip

\noindent
{\bf Remark} $\s\s$ Our algorithm also works for non-autonomous problems
\begin{eqnarray} \nonumber
& u''''(x)=\la f(x,u(x)), \s  \mbox{for $x \in (-1,1)$} \\ \nonumber
&  u(-1)=u'(-1)=u(1)=u'(1)=0 \,, \nonumber
\end{eqnarray}
although it is not known if the value of $u(0)$ is a global parameter here. In case $u(0)$ is not a global parameter, there might be solutions not lying on the computed solution curve.

\section{Numerical continuation  of solutions in the first harmonic}
\setcounter{equation}{0}
\setcounter{thm}{0}
\setcounter{lma}{0}

We describe numerical computation of solutions for the problem
\beq
\lbl{n1}
u''+f(u)=\mu \sin x+e(x), \s 0<x<\pi, \s u(0)=u(\pi)=0 \,.
\eeq
On the interval $(0,\pi)$, we have $\la _1=1$, $\p _1(x)=\sin x$, $\la _2=4$, $\p _2(x)=\sin 2x$, etc. Decompose $u(x)=\xi \sin x +U(x)$, with $\int _0^{\pi} U(x) \sin x \, dx=0$. It was shown in P. Korman \cite{K2014} that $\xi$ is a global parameter, uniquely identifying the solution pair $(u(x), \mu)$ of (\ref{sept18a}), provided that 
\beq
\lbl{n1aa}
f'(u)<\la _2=4 \,, \s \mbox{ for all $u \in R$} \,.
\eeq
In case the condition (\ref{n1aa}) fails, the value of $\xi$ is  not a global parameter in general, as we shall see later in this section.
However, we shall still compute the solution curve of (\ref{n1}): $(u(\xi),\mu (\xi))$,  although if the condition (\ref{n1aa}) fails, there might be solutions not lying on the computed solution curve.
We use   Newton's method to perform continuation in $\xi$.
\medskip

We begin by implementing the ``linear solver", i.e., the numerical solution of the following problem: given any $\xi \in R$, and any continuous functions $a(x)$ and $g(x)$, find $u(x)$ and $\mu \in R$ solving
\beqa
\lbl{n2}
& u''+a(x)u=\mu \sin x+g(x), \s 0<x<\pi \,, \\\nonumber
& u(0)=u(\pi)=0 \,,\\ \nonumber
& \int _0^{\pi} u(x) \sin x \, dx=\xi \,.\nonumber
\eeqa
Here $g(x)$ is any continuous function, not necessarily orthogonal to $\sin x$.
The general solution of the differential equation in (\ref{n2}) is
\beq
\lbl{n2a}
u(x)=Y(x)+c_1 u_1(x)+c_2 u_2(x) \,,
\eeq
where $Y(x)$ is any particular solution of that equation, and $u_1(x)$, $u_2(x)$ are two linearly independent solutions of the corresponding homogeneous equation
\beq
\lbl{n3}
u''+a(x)u=0, \s 0<x<\pi \,.
\eeq
We shall calculate the particular solution in the form  $Y(x)=\mu Y_1(x)+ Y_2(x)$, where $Y_1(x)$ solves 
\[
u''+a(x)u=\sin x, \s u(0)=0, \s u'(0)=1 \,,
\]
and   $Y_2(x)$ is the solution of
\[
u''+a(x)u=g(x), \s u(0)=0, \s u'(0)=1 \,.
\]
Let $u_1(x)$ be the solution of the equation (\ref{n3}) together with  the initial conditions $u(0)=0$, $u'(0)=1$, and let  $u_2(x)$ be any solution of (\ref{n3}) with  $u_2(0) \ne 0$.
The condition $u(0)=0$ in (\ref{n2}) implies that $c_2=0$ in (\ref{n2a}), so that  there is no need to compute $u_2(x)$, and (\ref{n2a}) takes the form  
\beq
\lbl{n4}
u(x)=\mu Y_1(x)+Y_2(x)+c_1 u_1(x) \,.
\eeq
The condition $u(\pi)=0$ and the last line in (\ref{n2}) imply that
\[
\mu Y_1(\pi)+c_1 u_1(\pi)=-Y_2(\pi) \,,
\]
\[
\mu \int _0^{\pi} Y_1(x)\sin x \, dx+c_1 \int _0^{\pi} u_1(x)\sin x \, dx=\xi -\int _0^{\pi} Y_2(x)\sin x \, dx \,,
\]
Solving this system for $\mu$ and $c_1$, and using $c_1$ in (\ref{n4}), gives the solution $(u(x),\mu)$ of (\ref{n2}).
\medskip

Turning to the problem (\ref{n1}), we begin with an initial value of $\xi =\xi _0$, and using a step size $\Delta \xi$, on a mesh $\xi _i=\xi _0 +i \Delta \xi$, $i=1,2, \ldots, nsteps$,  we compute the solution $(u(x),\mu)$ of (\ref{n1}), satisfying $\int _0^{\pi} u (x) \sin x \, dx=\xi _i$, by using Newton's method. Namely, assuming that the iterate $(u_n(x), \mu _n)$ is already computed, we linearize the equation (\ref{n1}) at $u_n(x)$:
\[
u''_{n+1}+f(u_n)+f'(u_n) \left(u_{n+1}-u_n \right)=\mu \sin x+g(x) \,,
\]
and calculate $(u_{n+1}(x),\mu _{n+1} )$ by using the method described above for 
 the problem (\ref{n2}) with
$a(x)=f'(u_n(x))$, $g(x)=-f(u_n(x))+f'(u_n(x)) u_n(x)+e(x)$, and $\xi =\xi _i$. After several iterations, we compute $(u(\xi _i), \mu (\xi _i))$. We found that three iterations of Newton's method, coupled with   $\Delta \xi$  not too large (say, $\Delta \xi=0.2$), were sufficient for accurate computation of the solution curves. To start Newton's iterations, we used the solution $u(x)$ computed at the preceding step, i.e., $u_0(x)=u(x)|_{\xi = \xi _{i-1}}$.
\medskip

\begin{figure}[h]
\begin{center}
\scalebox{0.65}{\includegraphics{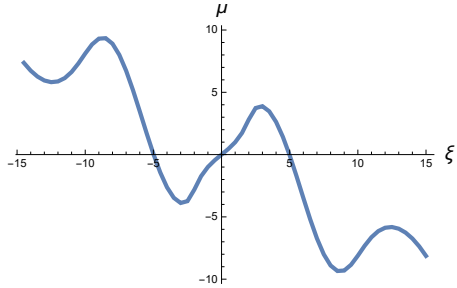}}
\end{center}
\caption{ The  global solution curve for  the  problem (\ref{n2new})}
\label{talk}
\end{figure}

\noindent
{\bf Example} $\s\s$ We used {\em Mathematica} to solve 
\beq
\lbl{n2new}
u''+\sin u=\mu \sin x+x-\frac{\pi}{2} \,, \s 0<x<\pi, \s u(0)=u(\pi)=0 \,.
\eeq
Observe that $\int_0^{\pi} \left(x-\frac{\pi}{2} \right) \sin x \, dx=0$. The solution curve $\mu =\mu (\xi)$ is presented in Figure \ref{talk}. The condition (\ref{n1aa}) applies here, and hence this curve exhausts the solution set of (\ref{n2new}).
\medskip

\medskip

\noindent
Let us look at the points of intersection of the solution curve with the $\xi$-axis, where $\mu=0$.  Figure \ref{talk}  strongly suggests that the problem
\[
u''+\sin u=x-\frac{\pi}{2} \,, \s 0<x<\pi, \s u(0)=u(\pi)=0
\]
has exactly three solutions, one of which has zero first harmonic.
\medskip

We now discuss computation of multiple solutions for the problem
\beq
\lbl{ca1}
u''+f(u)=0 \,, \s 0<x<\pi, \s u(0)=u(\pi)=0 \,.
\eeq
Of course, solutions of (\ref{ca1}) can be computed by shooting, but we describe a more general approach which is also applicable to PDE's.
Embed (\ref{ca1}) into a family of  problems
\beq
\lbl{ca2}
u''+f(u)=\mu _1 \sin x \,, \s 0<x<\pi, \s u(0)=u(\pi)=0 \,.
\eeq
Decompose its solution as $u(x)=\xi_1 \sin x +U(x)$, with $\int _0^{\pi} U(x) \sin x \, dx=0$. As above, we can compute the solution pair $(u(x),\mu _1)$ as a function of $\xi _1$, and draw the curve $\mu_1= \mu_1(\xi _1)$. At the points of intersection of this curve with the $\xi _1$ axis, where $\mu _1=0$, we obtain  solutions of (\ref{ca1}). Similarly, we can embed (\ref{ca1}) into another family of problems
\beq
\lbl{ca3}
u''+f(u)=\mu _2 \sin2  x \,, \s 0<x<\pi, \s u(0)=u(\pi)=0 \,.
\eeq
Decompose its solution as $u(x)=\xi_2 \sin2 x +U(x)$, with $\int _0^{\pi} U(x) \sin 2x \, dx=0$. We  perform the continuation in the second harmonic $\xi _2$ in a completely similar way to the computation in the first harmonic, described above. We compute the solution pair $(u(x),\mu _2)$ as a function of $\xi _2$, and draw the curve $\mu_2= \mu_2(\xi _2)$. At the points of intersection of this curve with the $\xi _2$ axis, where $\mu _2=0$, we obtain  solutions of (\ref{ca1}). Our computations indicate that these solutions of (\ref{ca1}) can be made  different from the ones obtained by the continuation in the first harmonic $\xi _1$, if $f'(0)> \la _2=4$.
\medskip

\begin{figure}[h]
\begin{center}
\scalebox{0.65}{\includegraphics{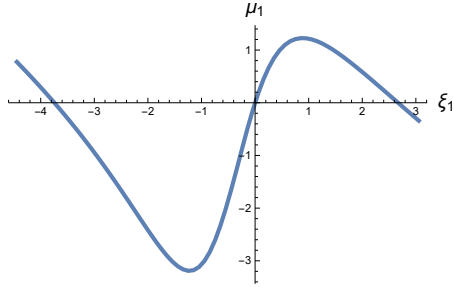}}
\end{center}
\caption{ The  solution curve $\mu_1= \mu_1(\xi _1)$ for  the  problem (\ref{ca5}) }
\label{castro1}
\end{figure}

In a recent paper A. Castro et al \cite{C} considered the problem (\ref{ca1}) (they also considered the corresponding PDE problem on a general domain). Assuming that $f(0)=0$, $f'(0)>\la _2$, $f'(0) \ne \la _k$ for any $k \geq 3$, and that there exists $\gamma<\la _1$ and $\rho >0$  such that $\frac{f(u)}{u} \leq \gamma$ for $|u| \geq \rho$, they proved that the problem (\ref{ca1}) has at least four non-trivial solutions, in addition to the trivial one. One of these solutions is positive, and one is negative. We now illustrate their result for the following example:
\beq
\lbl{ca4}
u''+\frac{6u}{1+u+2u^2}=0 \,, \s 0<x<\pi, \s u(0)=u(\pi)=0 \,.
\eeq
Here $f'(0)=6$ lies between $\la _2=4$ and $\la _3=9$. We calculated the solution curve $\mu_1= \mu_1(\xi _1)$ for the problem 
\beq
\lbl{ca5}
u''+\frac{6u}{1+u+2u^2}=\mu _1 \sin x \,, \s 0<x<\pi, \s u(0)=u(\pi)=0 \,.
\eeq
The result is presented in Figure \ref{castro1}.

\medskip

The function $\mu_1= \mu_1(\xi _1)$ has a root at $\xi _1 \approx -3.8$, which implies that the problem (\ref{ca4}) has a solution with the first harmonic 
$\xi _1 \approx -3.8$. This solution, which we denote by $u_1(x)$, is negative on $(0,\pi)$. We computed $u_1'(0) \approx -3.1968$, which together with $u_1(0)=0$ is sufficient to identify this solution. The root $\xi _1=0$ corresponds to the trivial solution of (\ref{ca4}). The root at $\xi _1 \approx 2.7$ gives us a second non-trivial solution  of (\ref{ca4}), $u_2(x)$, which is positive. From the solution $u_2(x)$ we find  $u_2'(0) \approx 2.0606$. 
\medskip

\begin{figure}[h]
\begin{center}
\scalebox{0.65}{\includegraphics{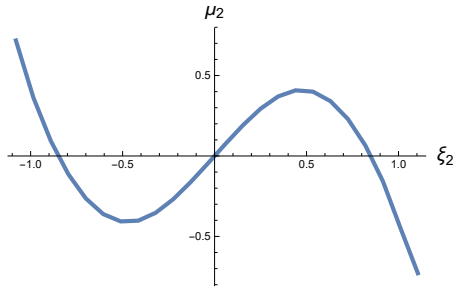}}
\end{center}
\caption{ The  solution curve $\mu_2= \mu_2(\xi _2)$ for  the  problem (\ref{ca5}) }
\label{k2curve}
\end{figure}

To obtain two more solutions of  (\ref{ca4}), we calculated the solution curve $\mu_2= \mu_2(\xi _2)$ for the problem 
\beq
\lbl{ca6}
u''+\frac{6u}{1+u+2u^2}=\mu _2 \sin 2x \,, \s 0<x<\pi, \s u(0)=u(\pi)=0 \,,
\eeq
as explained above.
The result is presented in Figure \ref{k2curve}. The function $\mu_2= \mu_2(\xi _2)$ has a root at $\xi _2 \approx -0.85$, which implies that the problem (\ref{ca4}) has a solution, called $u_3(x)$,  with the second harmonic $\xi _2 \approx -0.85$. This solution is sign-changing with exactly one root inside $(0,\pi)$, so that it is different from $u_1(x)$ and   $u_2(x)$. From the solution $u_3(x)$ we find  $u_3'(0) \approx -1.222$ The root $\xi _2=0$ again corresponds to the trivial solution of (\ref{ca4}). The root at $\xi _2 \approx 0.85$ gives us the fourth non-trivial solution  of (\ref{ca4}), $u_4(x)$. From the solution $u_4(x)$ we find  $u_4'(0) \approx 1.222$. In fact, $u_4(x)=u_3(\pi-x)$. (In general, if $u(x)$ is a solution of (\ref{ca4}), so is $u(\pi-x)$.) 

\medskip

When  doing continuation in $\xi _2$, one discovers that there is at least one more solution curve, in addition to the one in Figure \ref{k2curve}. It is the curve joining $u_1(x)$ and $u_2(x)$. To ``jump" on the solution curve in Figure \ref{k2curve}, we began Newton's iterations at each $\xi _2$ with $u_0(x)=\sin 2x$, a function that resembles the solutions $u_3(x)$ and $u_4(x)$.
\medskip

As we mentioned in the Introduction, $\xi _1$ is a global parameter for the problem (here $\int _0^{\pi} e(x) \sin x \, dx=0$)
\[
u''+f(u)=\mu _1 \sin x +e(x)\,, \s 0<x<\pi, \s u(0)=u(\pi)=0 \,,
\]
uniquely determining the solution pair $(\mu _1,u(x))$, provided that $f'(u)<\la _2=4$, for all $u \in R$. The problem  (\ref{ca5}) illustrates that this condition cannot be dropped, because  (\ref{ca5}) has at least two solution curves $\mu _1=\mu _1( \xi _1)$ (since  the solutions $u_3(x)$ and $u_4(x)$ do not lie on the solution curve in Figure \ref{castro1}).
 
\section{Some open problems}

Our computations raise a number of open questions, both of theoretical and computational nature (some  were already mentioned above).
\medskip

\begin{figure}[h]
\begin{center}
\scalebox{0.65}{\includegraphics{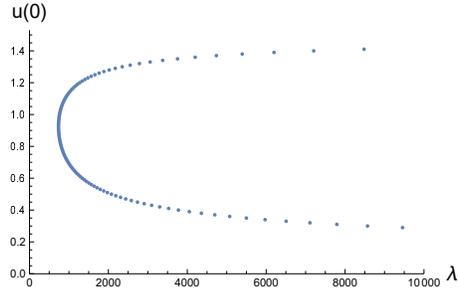}}
\end{center}
\caption{ The first solution curve  for  the  problem (\ref{ln1}) }
\label{firstcurve}
\end{figure}

\begin{figure}[h]
\begin{center}
\scalebox{0.65}{\includegraphics{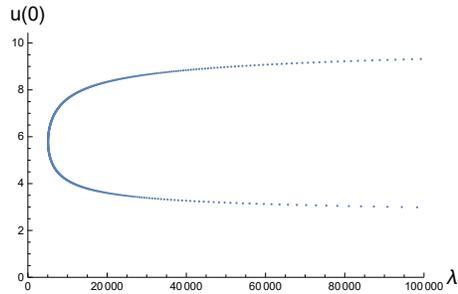}}
\end{center}
\caption{ The second solution curve  for  the  problem (\ref{ln1}) }
\label{secondcurve}
\end{figure}

\begin{figure}[h]
\begin{center}
\scalebox{0.65}{\includegraphics{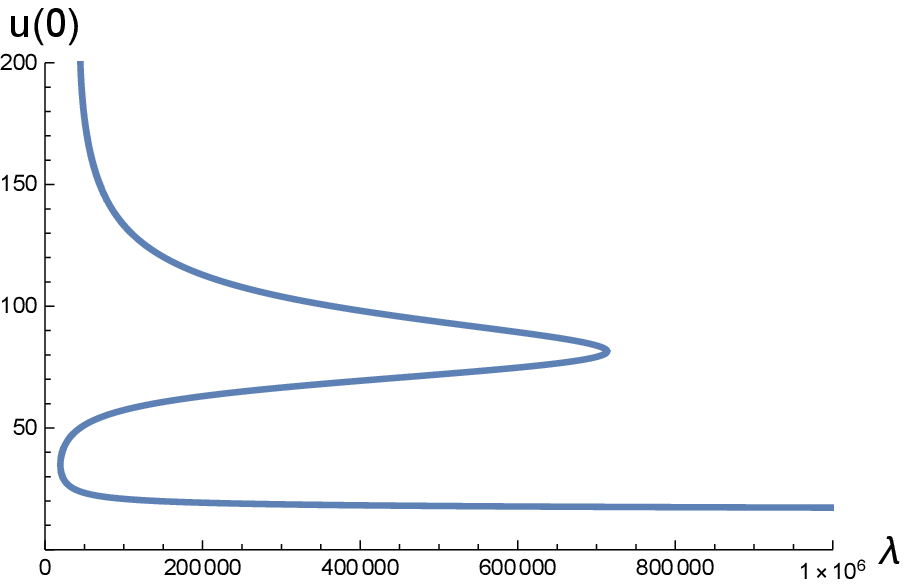}}
\end{center}
\caption{ The third (highest) solution curve  for  the  problem (\ref{ln1}) }
\label{thirdcurve}
\end{figure}

Equations with supercritical nonlinearities present both theoretical and computational challenges, which we explain on the example of the  equation due to C.S. Lin and W.-M. Ni \cite{LN} 
\beq
\lbl{ln1}
 u'' +\frac{n-1}{r}u'+\la \left( u^q+u^{2q-1} \right)=0 \,, \s  \; u'(0)=u(1)=0 \,,
\eeq
with $\frac{n}{n-2}<q<\frac{n+2}{n-2}<2q-1$. Here $\frac{n+2}{n-2}$ is the critical exponent for the space dimension $n$. When shooting, solutions of  supercritical equations may fail to achieve the root of either $u(r)$ or $u'(r)$ (leading to the {\em ground state solutions}), therefore computations were terminated only when the negative values of $u(r)$ were achieved.  Moreover, one can have several solution curves involving vastly different scales. 
We solved numerically the problem (\ref{ln1}) with $n=3$ and $q=4$. The problem has two parabola-like curves opening to the right in the $(\la,u(0))$ plane, one above the other, given in Figure \ref{firstcurve} and Figure \ref{secondcurve}, and an even higher third curve given in Figure \ref{thirdcurve}, at much larger values of $\la$. All three curves exhibit horizontal asymptotes. The horizontal asymptotes are not possible for subcritical $f(u)$, with $f(u)>0$, see T. Ouyang and J. Shi \cite{OS}. A very challenging open problem is to show that there is exactly one turn on the lower two curves of (\ref{ln1}). 

\medskip

When   computing  the curves of positive solutions of non-autonomous problems
\[
 u'' +\frac{n-1}{r}u'+\la f(r,u) =0 \,, \s  \; u'(0)=u(1)=0 \,,
\] 
it is desirable to know under what conditions  the value of $u(0)$ is a global parameter, as we saw above. In case $n=1$, $u(0)$ is known to be  a global parameter, assuming that the condition $f_r (r,u)\leq 0$ holds, see P. Korman \cite{K2013}, or P. Korman and J. Shi \cite{KS1}. We conjecture that the same result holds in case $n>1$ too. It would be desirable to find other conditions for this result to hold, and some counter-examples.
\medskip

Similarly, for positive solutions of the problem
\[
 u'' +\frac{n-1}{r}u'+\la f(u)+g(u) =0 \,, \s  \; u'(0)=u(1)=0  
\]
it would be desirable to find conditions guaranteeing  that $u(0)$ is  a global parameter. Then solution curves can be computed by a method similar to the one we used for non-autonomous problems. In case $f(u)$ and $g(u)$ are pure powers, one can rescale the problem, putting the parameter in front of the nonlinearity, and apply the shoot-and-scale method.
\medskip

For  positive solutions  of the elastic beam equation,  clamped at the end points,
\begin{eqnarray} \nonumber
& u''''(x)=\la f(u(x)), \s  \mbox{for $x \in (-1,1)$} \\ \nonumber
&  u(-1)=u'(-1)=u(1)=u'(1)=0 \,, \nonumber
\end{eqnarray}
it  would be interesting to study further under what conditions  the value of $u(0)$, the maximal value of solution, is a global parameter. Can one drop the conditions (\ref{bm2})? Are there any counterexamples? It appears to be  a very challenging research direction to prove the exact multiplicity results suggested by our computations. For example, is the solution curve for the problem (\ref{bm20}) exactly $S$-shaped? The time map method does not apply here, and the bifurcation approach is not sufficiently developed.
\medskip

Another case of  a global parameter occurs for positive solutions $u=u(x,y)$ of  the problem
\[
u_{xx}+u_{yy}+\la f(u)=0 \s \mbox{for $(x,y) \in D$}  \,, \; u=0 \; \mbox{on $\partial D$} \,,
\]
where $D$ is a class of symmetric domains in the $(x,y)$-plane, which in particular includes ellipses (see \cite{HK}, or \cite{K2}).
By the results of B. Gidas, W.-M. Ni and L. Nirenberg \cite{GNN} any   positive solution is symmetric with respect to $x$ and $y$, so that $u(0,0)$ gives the maximum 
value of any positive solution. M. Holzmann   and  H. Kielhofer \cite{HK} showed that  $u(0,0)$  is in fact a global parameter, uniquely identifying the solution pair $(\la ,u(x,y))$. It would be interesting to compute the global solution curves in the $(\la, u(0,0))$-plane.
\medskip

In a recent paper \cite{KS2}, we provided both theoretical and computational results for first order equations, periodic in time, with the approach similar to that of Section 8. A number of open questions were raised in that paper.


\begin{thebibliography}{99}
\bibitem{Al0}
E.L. Allgower,   On a discretization of $y \, ''+\lambda y\sp{k}=0$, {\em Topics in numerical analysis}, {\rm II} (Proc. Roy. Irish Acad. Conf., Univ. College, Dublin, 1974), pp. 1-15. (Academic Press), London  (1975).
\vspace{-0.2cm}  

\bibitem{A}
E.L. Allgower  and K. Georg,   Numerical Continuation Methods. An Introduction. Springer Series in Computational Mathematics, {\bf 13}. Springer-Verlag, Berlin  (1990).
\vspace{-0.2cm} 

\bibitem{C}
A. Castro, J. Cossio, S. Herr\'{o}n and C. V\'{e}lez,  Existence and multiplicity results for a semilinear elliptic problem, {\em  J. Math. Anal. Appl.} {\bf  475}, no. 2, 1493-1501  (2019).
\medskip

\bibitem{D}
 E.N. Dancer, On the structure of solutions of an equation in catalysis theory when a parameter is large, {\em   J. Differential Equations} {\bf  37}, no. 3,  404-437    (1980). 
 \vspace{-0.2cm}

\bibitem{GNN}
B. Gidas, W.-M. Ni and L. Nirenberg,   Symmetry and related properties via the      
     maximum principle, {\em Commun. Math. Phys}. {\bf 68}, 209-243 (1979).
\vspace{-0.2cm}

\bibitem{HK}
M. Holzmann   and  H. Kielhofer, Uniqueness of global positive solution branches of nonlinear elliptic problems, {\em Math. Ann.} {\bf 300},  no. 2,  221-241   (1994).
\vspace{-0.2cm}

\bibitem{I}
J.A. Iaia,   Localized solutions of elliptic equations: loitering at the hilltop, {\em  Electron. J. Qual. Theory Differ. Equ.}, No. 12  (2006). 
\vspace{-0.2cm}

\bibitem{JL}
 D.D. Joseph and    T.S. Lundgren,
Quasilinear Dirichlet problems driven by positive sources, 
{\em Arch. Rational Mech. Anal.} {\bf 49},   241-269 (1972/1973).
\vspace{-0.2cm}

\bibitem{K}
P. Korman,  Uniqueness and exact multiplicity of solutions for a class of fourth-order semilinear problems, {\em  Proc. Roy. Soc. Edinburgh Sect. A} {\bf 134}, no. 1, 179-190  (2004).
\vspace{-0.2cm}

\bibitem{K4}
P. Korman,  Computation of radial solutions of semilinear equations, {\em   Electron. J. Qual. Theory Differ. Equ.}, No. 13  (2007).
\vspace{-0.2cm}

\bibitem{K17}
P. Korman,  An oscillatory bifurcation from infinity, and from zero, {\em   NoDEA Nonlinear Differential Equations Appl.} {\bf   15},  no. 3, 335-345    (2008).
\vspace{-0.2cm}

\bibitem{K2}
P. Korman, Global Solution Curves for Semilinear Elliptic Equations, World Scientific, Hackensack, NJ (2012).
 \vspace{-0.2cm}

\bibitem{K2013}
P. Korman,   Exact multiplicity and numerical computation of solutions for two classes of non-autonomous problems with concave-convex nonlinearities, {\em Nonlinear Anal.} {\bf 93}, 226-235 (2013).
\vspace{-0.2cm}

\bibitem{K2014}
P. Korman,  Curves of equiharmonic solutions, and problems at resonance, {\em Discrete Contin. Dyn. Syst.} {\bf 34}, no. 7, 2847-2860  (2014).
\vspace{-0.2cm}

\bibitem{K2015}
P. Korman,  Regularization of radial solutions of $p\, $-Laplace equations, and computations using infinite series, {\em  Electron. J. Qual. Theory Differ. Equ. } No. 40 (2015).
 \vspace{-0.2cm}

\bibitem{KLS} P.  Korman, Y. Li  and  D.S. Schmidt,  A computer assisted study of uniqueness of nodal ground state solutions, {\em  J. Comput. Appl. Math.} {\bf  356}, 402-406  (2019). 
 \vspace{-0.2cm}

\bibitem{KS2}
P.  Korman  and  D.S. Schmidt, Global solution curves for first order periodic problems, with applications. Preprint.
\vspace{-0.2cm}

\bibitem{KS}
P. Korman and  J. Shi,   On Lane-Emden type systems, {\em  Discrete Contin. Dyn. Syst.},  suppl.,  510-517 (2005).
\vspace{-0.2cm}

\bibitem{KS1}
P. Korman  and  J. Shi,   Instability and exact multiplicity of solutions of semilinear equations, Proceedings of the Conference on Nonlinear Differential Equations (Coral Gables, FL, 1999), 311-322 (electronic), Electron. J. Differ. Equ. Conf., 5, Southwest Texas State Univ., San Marcos, TX  (2000).
 \vspace{-0.2cm}

\bibitem{LN}
C.S. Lin  and  W.-M. Ni  A counterexample to the nodal domain conjecture and related semilinear equation, {\em Proc. Amer. Math. Soc.} {\bf 102},  271-277   (1988).
\vspace{-0.2cm}

\bibitem{OS}
T. Ouyang and J. Shi,
Exact multiplicity of positive solutions for a class of semilinear problems, II, {\em J. Differential Equations} {\bf 158},  no. 1, 94-151 (1999).
\vspace{-0.2cm}

\bibitem{PS1}
L.A. Peletier and  J. Serrin, Uniqueness of positive solutions of semilinear equations in $ R\sp{n}$, {\em Arch. Rational Mech. Anal.} {\bf 81}  no. 2,  181-197   (1983).
\vspace{-0.2cm} 

\bibitem{QS}
P. Quittner  and   P. Souplet,  Superlinear Parabolic Problems. Blow-up, global existence and steady states. Birkh\"{a}user Advanced Texts: Basler Lehrb\"{u}cher.  (Birkh\"{a}user), Basel   (2007).
\vspace{-0.2cm}

\bibitem{S10}
R. Schaaf,    Global Solution Branches of Two Point Boundary Value Problems, Lecture Notes in 
Mathematics, no. {\bf 1458}, (Springer-Verlag) (1990).
\vspace{-0.2cm}

\end{thebibliography}
\end{document}